\documentclass[12pt]{article}
\usepackage{amssymb,latexsym,theorem}

\newcommand{\U}{{\cal U}}
\newcommand{\cl}{{\cal L}}
\newcommand{\cm}{{\cal M}}
\newcommand{\cs}{{\cal S}}
\newcommand{\suppd}{{\mathrm {supp}_\nabla}}

\newcommand{\res}{\mathop{\mathrm{res}}\nolimits}
\newcommand{\ind}{\mathop{\mathrm{ind}}\nolimits}
\newcommand{\St}{\mathop{\mathrm{St}}\nolimits}
\newcommand{\End}{\mathop{\mathrm{End}}\nolimits}
\newcommand{\Ext}{\mathop{\mathrm{Ext}}\nolimits}
\newcommand{\sEnd}{\mathop{{\cal E}nd}\nolimits}
\newcommand{\EndF}{\End_F}
\newcommand{\sEndF}{\sEnd_F}

\newcommand{\Es}{$\mathrm{E}_6$}
\newcommand{\Ff}{$\mathrm{F\!}_4$}

\newcommand{\qed}{\unskip\nobreak\hfill\hbox{ $\Box$}}

\newtheorem{Proposition}[subsection]{Proposition}
\newtheorem{Theorem}[subsection]{Theorem}
\newtheorem{Lemma}[subsection]{Lemma}
\newtheorem{Corollary}[subsection]{Corollary}
\theorembodyfont{\normalfont}
\newtheorem{Remark}[subsection]{Remark}
\newtheorem{Example}[subsection]{Example}
\newtheorem{Definition}[subsection]{Definition}

\begin{document}
\title{Steinberg modules and Donkin pairs.}
\author{Wilberd van der Kallen}
\maketitle
\sloppy

\section{Summary}
We prove that in characteristic $p>0$ a module with good filtration
for a group of type {\Es} restricts to a module with good filtration for
a group of type {\Ff}. Thus we confirm a conjecture of Brundan for one more
case. Our method relies on the canonical Frobenius splittings of Mathieu.
Next we settle the remaining cases, in characteristic not 2,
with a computer-aided variation on the old method of Donkin.

\section{Preliminaries}
Our base field $k$ is algebraically closed of characteristic $p$.
Let $G$ be a connected semisimple group  and $H$ a connected
semisimple subgroup. (Good filtrations with more general groups are treated
in \cite{Donkin}.) We refer to \cite{Jantzen}
and \cite{vdK}
for unexplained terminology and notation.

Now choose a Borel subgroup $B$ in $G$ and a maximal torus $T$ in $B$ so
that, if $B^-$ is the opposite Borel subgroup, then
$B\cap H$ and $B^-\cap H$ are a Borel subgroups in $H$ and
% thus?
$T\cap H$ is a maximal torus in $H$.
% At least the intersection of two Borel subgroups contains a maximal torus.
%
% Choose the tori first. First a maximal torus $S$ in $H$, then
% (in its centralizer in $G$) the maximal torus $T$. The intersection of $T$
% with $H$ can not be too large, as $S$ is its own centralizer in $H$.
% Then choose a generic cocharacter of $T\cap H=S$.
% It defines an order on both  character groups. We get a Borel group $B_H$
% in $H$ and a parabolic $P$ in $G$ so that $B_H$ is in the unipotent
% radical of $P$ and its opposite is in the unipotent radical of the
% opposite parabolic. The Levi subgroup is the centralizer we mentioned
% earlier.

We follow the convention that the roots of $B$ are positive.
If $\lambda\in X(T)$ is dominant, then $\ind_B^G(-\lambda)$ is the
dual Weyl module
$\nabla_G(\lambda^*)$ with highest weight $\lambda^*=-w_0\lambda$ and lowest
weight
$-\lambda$. Its dual is the Weyl module $\Delta_G(\lambda)$.
In a good filtration of a $G$-module the layers are of the form
$\nabla_G(\mu)$.

\begin{Definition}
We say that $(G,H)$
is a  Donkin pair if for any $G$-module $M$ with good filtration,
the $H$-module
$\res^G_HM$ has good filtration.
\end{Definition}

Let $\U(U)$ denote the hyperalgebra of the unipotent radical $U$ of $B$.
We recall the presentation of Weyl modules.

\begin{Lemma}\label{Weylpres}
Let $\lambda$ be dominant and let
$v_{-\lambda^*}$ be a nonzero weight vector of lowest weight $-\lambda^*$ in
$\Delta_G(\lambda)$. Then  $v_{-\lambda^*}$ generates $\Delta_G(\lambda)$
as a $\U(U)$-module, and the annihilator of $v_{-\lambda^*}$ equals the
left ideal of $\U(U)$ generated by the $X_\alpha^{(n)}$
with $\alpha$ simple and $n>(\lambda^*,\alpha^\vee)$.
\end{Lemma}

\paragraph{Proof}
Note that $\U(U)$ is a graded algebra graded by height.
Therefore the left ideal
% $F$
in the lemma is the intersection of all ideals
$I$ of finite codimension that contain it and that lie inside
the annihilator.
% Take $I=F+\{x\mid ht(x)>M\}$
But by the proof of \cite[Proposition Fondamentale]{Polo var}
such ideals $I$ are equal to the annihilator.
\qed\vskip\baselineskip

Let $X$ be a smooth projective $B$-variety with canonical bundle $\omega$.
(Generalizations to other varieties will be left to the reader.)
There is by \cite[\S 2]{Mehta-Ramanathan}
a natural map $\epsilon:H^0(X,\omega^{1-p})\to k$ so that
$\phi\in H^0(X,\omega^{1-p})$ determines a Frobenius
splitting if and only if
$\epsilon(\phi)=1$.
Let $\St_G$ be the Steinberg module of the simply connected cover $\tilde G$
of
$G$.
For simplicity of notation we further assume that $\St_G$ is
actually a $G$-module. Its $B$-socle is the
highest weight space $k_{(p-1)\rho}$.

Recall that a Frobenius splitting of $X$
is called \emph{canonical} if the corresponding $\phi$ is $T$-invariant and
lies in the image of a $B$-module map
$\St_G\otimes k_{(p-1)\rho}\to H^0(X,\omega^{1-p})$.
(Compare lemma \ref{Weylpres} and \cite[Definition 4.3.5]{vdK}.)
If the group $G$ needs to be emphasized, we will speak of a $G$-canonical
splitting.
Now suppose $X$ is actually a $G$-variety.

\begin{Lemma}\label{StSt}
$X$ has a canonical splitting if and only if there is a $G$-module
map $\psi:\St_G\otimes\St_G\to H^0(X,\omega^{1-p})$ so that
$\epsilon\psi\neq0$.
\end{Lemma}

\paragraph{Proof}
There is, up to scalar multiple, only one possibility for a map
$\St_G\otimes\St_G\to k$.
If $\epsilon\psi\neq0$ then the subspace of $T$-invariants in
$\St_G\otimes k_{(p-1)\rho}$ maps isomorphically to $k$. Conversely, a map
from $\St_G\otimes k_{(p-1)\rho}$ to a $G$-module $M$ can be extended
to $\St_G\otimes\St_G$ because the $G$-module generated by the image
of $k_{(p-1)\rho}$ in $M\otimes\St_G^*$ is $\St_G$.
\qed\vskip\baselineskip

We have the following fundamental result of Mathieu \cite{Mathieu G}.
\begin{Theorem}\textup{\cite[6.2]{Mathieu T}}\label{splitgood}
Assume $X$ has a canonical splitting and $\cl$ is a $G$-linearized
line bundle on $X$. Then $H^0(X,\cl)$ has a good filtration.
\end{Theorem}

\section{Pairings}
Now apply this to $X=G/B$. Of course the $\nabla_G(\mu)$
are of the form $H^0(X,\cl)$, see \cite[I 5.12]{Jantzen}.
It follows that $(G,H)$ is a  Donkin pair if
$X$ has an $H$-canonical splitting.
We also have a surjection $\St_G\otimes\St_G\to H^0(X,\omega^{1-p})$,
by \cite[II 14.20]{Jantzen}.
The composite with $H^0(X,\omega^{1-p})\to k$
may be identified as in \cite{Lauritzen Thomsen}, \cite{Mehta-Venkataramana}
 with the natural pairing on the self dual representation
$\St_G$. Thus we get

\begin{Theorem}[Pairing criterion]\label{pairing criterion}
Assume there is an $H$-module map
$$\St_H^*\otimes \St_H\to \St_G^*\otimes\St_G$$
whose composite with the evaluation map $\St_G^*\otimes\St_G\to k$
is nonzero.
Then $(G,H)$ is a Donkin pair.
\end{Theorem}

\begin{Remark}
Despite the notation, the $\tilde H$-module
$\St_H$ need not be an $H$-module.
Even if $\St_H$ is not an $H$-module, $\St_H^*\otimes \St_H$ is one.
It may be better to replace $H\subset G$ with the homomorphism
$\tilde H\to \tilde G$. Thus an operation like $\res_H^G$ would really mean
restriction along $\tilde H\to \tilde G$.
\end{Remark}

\begin{Remark}
The pairing criterion
is satisfied if and only if $G/B$ has an $H$-canonical splitting.
Indeed suppose we are given a map
$\St_H^*\otimes\St_H=\St_H\otimes\St_H\to H^0(G/B,\omega^{1-p})$
as in Lemma \ref{StSt}.
We have to factor it through
the surjection $\pi:\St_G^*\otimes\St_G\to H^0(G/B,\omega^{1-p})$.
But the kernel $K$ of $\pi$ has good filtration by
\cite{Mathieu G} (or by the proof in \cite[II 4.16]{Jantzen}),
so $\Ext_H^1(\St_H^*\otimes\St_H,\res_H^G K)$ vanishes by
theorem \ref{splitgood} and the main properties
of good filtrations (\cite[Theorem 1]{Mathieu G}, \cite[II 4.13]{Jantzen}).
\end{Remark}

Now we illustrate the criterion with some old examples of  Donkin pairs.

\begin{Example}\label{products}
Let $G$ still be semisimple and connected.
It is easy to see from the formulas in the proof of \cite[3.2]{Lauritzen
Thomsen}
that the pairing criterion applies to the diagonal $G$ inside a product
$G\times\cdots\times G$. \end{Example}

\begin{Example}\label{levi}
Let $H$ be the commutator subgroup of a Levi subgroup of a parabolic in
the semisimple connected
group $G$.
Then, after passing to simply connected covers if necessary,
% it is not necessary if one gives strange meaning to $\res^G_H$
% as in earlier remark.
$\St_H$ is a direct summand of $\res^G_H\St_G$, so again the
pairing criterion applies.
\end{Example}

\begin{Lemma}
Let $(G,H)$ satisfy the pairing criterion and let $X$ be a smooth projective
$G$-variety. If $X$ has a $G$-canonical splitting, then it has an
$H$-canonical one.
\end{Lemma}
\paragraph{Proof} Use lemma \ref{StSt}.\qed\vskip\baselineskip

The following lemma was pointed out to me by Jesper Funch Thomsen.

\begin{Lemma}\label{Gproduct}
Let  $X$, $Y$ be smooth projective
$G$-varieties with canonical splitting. Then $X\times Y$ has a $G$-canonical
splitting.
\end{Lemma}

\paragraph{Proof} Use example \ref{products}.\qed

\begin{Remark}
For the users of our book, let us now point out how to get theorem
\ref{splitgood}.
We have $G\times^BX=G/B\times X$ by remark \cite[1.2.2]{vdK},
so \cite[lemma 4.4.2]{vdK} applies with $Y=X$ in the notations of
that lemma.
\end{Remark}

\begin{Remark}
In lemma \ref{Gproduct} one cannot replace $G$ with $B$.
Here is an example. Take $G=SL_3$ in characteristic 2
and let $Z$ be the Demazure resolution
of a Schubert divisor. Then $H^0(Z,\omega^{-1}_Z)$ is a nine dimensional
$B$-module. There is a fundamental representation $V$ so that
$H^0(Z,\omega^{-1}_Z)$ is isomorphic to a codimension one
submodule of the
degree three part of the ring
of regular functions on $V$.
Using this,
one checks with computer assisted computations that $Z$, $Z\times Z$,
$Z\times Z\times Z$ have $B$-canonical splittings, while
$Z\times Z\times Z\times Z$ does not have one.
\end{Remark}

Our next aim is to treat the following example.
\begin{Example}\label{E6-F4}
For $G$ we take the simply connected group of type $E_6$. From the symmetry
of its Dynkin diagram we have a graph automorphism which is an involution.
For $H$ we take the group
of fixed points of the involution.
It is connected (\cite[8.2]{Steinberg}) of type $F_4$.
It has been conjectured by Brundan \cite[4.4]{Brundan} that $(G,H)$ is a
Donkin pair.
\end{Example}

More generally, with our usual notations we have.

\begin{Theorem}\label{three pieces}
Assume there are
dominant weights $\sigma_1$, $\sigma_2$, $\sigma_3$, so that
\begin{enumerate}
\item The highest weight $(p-1)\rho_G$
of $\St_G$ equals
$\sigma_1+\sigma_2+\sigma_3$.
\item $\sigma_1+\sigma_2$ and $\sigma_2+\sigma_3$ both
restrict to the highest weight $(p-1)\rho_H$ of $\St_H$.
\item The natural map
$\nabla_G(\sigma_1)\to\nabla_H(\res^{B}_{B\cap H}\sigma_1)$
is surjective.
\end{enumerate}
Then $(G,H)$ is a  Donkin pair.
In fact it satisfies the pairing criterion.
\end{Theorem}

\begin{Remark}
If $(G,H)$ is a  Donkin pair and $\lambda$ is dominant,
then one knows that
$\nabla_G(\lambda)\to\nabla_H(\res^{B}_{B\cap H}\lambda)=
\ind_{H\cap B^-}^H(\res^{B^-}_{H\cap B^-}\lambda)$,
induced by the projection of $\nabla_G(\lambda)$ onto its highest weight space,
is surjective. (Exercise. Use a good filtration as in
the proof of \cite[II 4.16]{Jantzen}.)
% By definition the highest weight vector survives.
% Factor out the earlier layers. We do not care if there are later layers.
% I do not quite understand what happens to lengths under restriction
% when you are looking at a more general $G$, $H$ pair.
% But I have made sure now that the projection onto the highest weight space
% of $\nabla_G(\lambda)$ is equivariant for the opposite Borel subgroup of $H$.
% For example, take $G$ of type $A_2$ and $H$ the $SL_2$ corresponding
% with a simple root. Now look what it says if $\lambda$ is the highest root.
\end{Remark}

\begin{Remark}
Our theorem \ref{three pieces} also applies to the Levi subgroup case of
example \ref{levi} (take $\sigma_1=0$).
One hopes to find a more general method to attack at least
all graph automorphisms. Theorem \ref{three pieces}
applies if the graph automorphism is an involution
and different simple roots in an orbit are perpendicular
to each other.
But for the graph automorphism of a group of type $\mathrm{A}_{2n}$
in characteristic $p>2$ there
are no $\sigma_1$, $\sigma_2$, $\sigma_3$ as in the theorem.
The coefficient of $\res^{B}_{B\cap H}\rho_G$ with respect to the fundamental
weight that corresponds to the short root is four, which is too high.
\end{Remark}

\subsubsection*{Proof of Theorem \ref{three pieces}.}
We will often write the restriction of a weight to $T\cap H$ with
the same symbol as the weight. We will repeatedly use basic properties
of Weyl modules and their duals.
See \cite[II 14.20]{Jantzen} for surjectivity of cup product between
dual Weyl modules
and \cite[II 2.13]{Jantzen} for Weyl modules as universal highest weight
modules.
We first need a number of nonzero maps of $H$-modules.
They are natural up to nonzero scalars that do not interest us.

The first map is the map $$\epsilon_H:\nabla_H(2(p-1)\rho_H)\to k$$
which detects Frobenius splittings on $H/(H\cap B)$. Together with the
surjection $$\nabla_H(\sigma_2)\otimes \nabla_H((p-1)\rho_G)\to
\nabla_H(2(p-1)\rho_H)$$ it gives a nonzero map
$\nabla_H(\sigma_2)\otimes \nabla_H((p-1)\rho_G)\to k$ and hence
a nonzero $$\eta_1:\nabla_H(\sigma_2)\to
\nabla_H((p-1)\rho_G)^*.$$

The map $\nabla_G(\sigma_2+\sigma_3)\to \St_H$ is nonzero, hence surjective.
The map $\nabla_G(\sigma_1)\to \nabla_H(\sigma_1)$ is surjective
by assumption.
In the commutative diagram
$$ \begin{array}{ccc}
 {\nabla_G(\sigma_2+\sigma_3)\otimes\nabla_G(\sigma_1)} & {\longrightarrow} &
\St_G \\[.2em]
 {\downarrow} & & {\downarrow}\\
 {\St_H\otimes\nabla_H(\sigma_1)} & {\longrightarrow} & {\nabla_H((p-1)\rho_G)}
\end{array}
$$
the horizontal maps are also surjective. So the map
$$\eta_2:\nabla_H((p-1)\rho_G)^*\to \St_G^*$$ is injective.
We obtain a nonzero $$\eta_2\eta_1:\nabla_H(\sigma_2)\to\St_G^*.$$
The nonzero $\St_H\to \nabla_G(\sigma_2+\sigma_3)$
combines with the map $$\nabla_G(\sigma_1)\otimes\nabla_G(\sigma_2+\sigma_3)
\to \St_G$$ to yield
$$\nabla_G(\sigma_1)\otimes \St_H\to \St_G$$ and combining this with
$\eta_2\eta_1$
we get
$$\eta_3: \nabla_H(\sigma_2)\otimes \nabla_G(\sigma_1)\otimes \St_H
\to \St_G^*\otimes\St_G.$$
We claim that its image is detected by the evaluation map
$$\eta_4:\St_G^*\otimes\St_G\to k.$$
This is because $\eta_3$ factors through $\nabla_H((p-1)\rho_G)^*\otimes
\St_G$, on which the restriction of $\eta_4$ factors through
$\nabla_H((p-1)\rho_G)^*\otimes
\nabla_H((p-1)\rho_G)$,
the map $\eta_1$ is nonzero, the image of
$ \nabla_G(\sigma_1)\otimes \St_H\to \St_G$ maps onto $\nabla_H((p-1)\rho_G)$.

{}From the nontrivial $\eta_4\eta_3$ we get a nontrivial
$$\eta_5:\nabla_H(\sigma_2)\otimes \nabla_G(\sigma_1)\to \St_H^*.$$
Then $\eta_5$ must be split surjective. Choose a left inverse
$$\eta_6:\St_H^*\to \nabla_H(\sigma_2)\otimes \nabla_G(\sigma_1)$$
of $\eta_5$.
It leads to
$$\eta_7:\St_H^*\otimes\St_H\to
\nabla_H(\sigma_2)\otimes \nabla_G(\sigma_1)\otimes\St_H$$
and the map we use in the pairing criterion is $\eta_3\eta_7$.
Indeed the map $\St_H^*\to\St_H^*$ defined by $\eta_4\eta_3\eta_7$
equals $\eta_5\eta_6$, hence is nonzero.
\qed\par

\section{The {\Es}-{\Ff} pair.}
We turn to the {\Es}-{\Ff} pair of example \ref{E6-F4}.
First observe that for $p>13$ one could simply follow the method of
\cite{Brundan} to prove that the pair is a Donkin pair.
Indeed the restriction to {\Ff} of a fundamental representation then
has its dominant weights in the bottom alcove.
Looking a little closer and applying the linkage principle one can treat
$p\geq 11$ in the same manner.
% For $p=11$ and $p=13$ the modules  $\nabla_H(\lambda)$ are irreducible
%  for
% $\lambda$ in $\{0,\varpi_1,\varpi_2,\varpi_3,
% \varpi_4,\varpi_1+\varpi_4\}$ (Bourbaki notation for {\Ff}),
% because those dominant weights that are distinct from $\lambda$
% lie in a different orbit
% than $\lambda$. It follows that the restriction
% to {\Ff} of a fundamental representation of {\Es} has a good filtration.

But for $p=5$ one has $\varpi_4\uparrow\varpi_1+\varpi_4$
and for $p=7$ one has $\varpi_1\uparrow\varpi_1+\varpi_4$.
% Let us write $(a,b,c,d)$ for $-\rho+
% a\epsilon_1+b\epsilon_2+c\epsilon_3+d\epsilon_4$.
% Then $0=(11/2,5/2,3/2,1/2)$,
% $-\rho+\varpi_1=(1,1,0,0)$,
% $-\rho+\varpi_2=(2,1,1)$,
% $-\rho+\varpi_3=(3/2,1/2,1/2,1/2)$,
% $-\rho+\varpi_4=(1,0,0,0)$.
% The fundamental alcove is given by $2a<p$, $a\geq b+c+d$,
% $b\geq c\geq d\geq0$.
% For $p=5$ one has $\varpi_4=(13/2,5/2,3/2,1/2)\uparrow
% (15/2,7/2,3/2,1/2) =\varpi_1+\varpi_4$.
% For $p=7$ one has $\varpi_1=(13/2,7/2,3/2,1/2)\uparrow
% (15/2,7/2,3/2,1/2) =\varpi_1+\varpi_4$.
This makes that one has more trouble to see that the restriction of
$\nabla_G(\varpi_4)$ has a good filtration with respective layers
$\nabla_H(\varpi_1)$, $\nabla_H(\varpi_3)$, $\nabla_H(\varpi_3)$,
$\nabla_H(\varpi_1+\varpi_4)$, $\nabla_H(\varpi_2)$.
% \begin{verbatim}branch([0,0,0,1,0,0],F4,*res_mat(F4,E6),E6)
%      2X[0,0,1,0] +1X[0,1,0,0] +1X[1,0,0,0] +1X[1,0,0,1]
% \end{verbatim}
For $p=2$ or $p=3$ it is even worse.

So let us apply theorem \ref{three pieces} instead.
We take $\sigma_1=(p-1)(\varpi_1+\varpi_3)$,
$\sigma_2=(p-1)(\varpi_2+\varpi_4)$, $\sigma_3=(p-1)(\varpi_5+\varpi_6)$
in the notations of Bourbaki for {\Es} \cite[Planches]{Bourbaki}.
Then $\res^{B}_{B\cap H}\varpi_i$ equals $\varpi_4$,
$\varpi_1$,
$\varpi_3$,
$\varpi_2$,
$\varpi_3$,
$\varpi_4$ for $i=1, \ldots ,6$
respectively.

First let let $p=2$. Then
$\nabla_H(\res^{B}_{B\cap H}\sigma_1)=\nabla_H(\varpi_3+\varpi_4)$
is irreducible. Indeed its dominant weights come in two parts.
% \begin{verbatim}dom_char(X[0,0,1,1],F4)
%    64X[0,0,0,0] +40X[0,0,0,1] + 8X[0,0,0,2] +14X[0,0,1,0] +
%     1X[0,0,1,1] + 2X[0,1,0,0] +24X[1,0,0,0] + 4X[1,0,0,1]
% \end{verbatim}
The weights $0$, $\varpi_4$, $\varpi_1$, $\varpi_3$, $2\varpi_4$,
 $\varpi_1+\varpi_4$, $\varpi_2$
lie in one orbit, and the highest weight lies in a different orbit
under the affine Weyl group.
To be more specific,
$\varpi_1-\rho_H\uparrow\varpi_3+\varpi_4$, but
$\varpi_4-\rho_H\uparrow0\uparrow\varpi_4\uparrow\varpi_1\uparrow\varpi_3
\uparrow2\varpi_4\uparrow\varpi_1+\varpi_4\uparrow\varpi_2$.
So $\nabla_G(\sigma_1)\to\nabla_H(\res^{B}_{B\cap H}\sigma_1)$
is surjective.
% $\varpi_3=(11/2,5/2,3/2,1/2)+(3/2,1/2,1/2,1/2)=
% (7,3,2,1)$,
% $\varpi_1-\rho=(1,1,0,0)\uparrow(7,3,2,0)
% \uparrow (8,3,2,1)=\varpi_3+\varpi_4$,
% $\varpi_4-\rho=(1,0,0,0)
% \uparrow(5,2,0,0)\uparrow(11/2,5/2,1/2,1/2)
% \uparrow(11/2,5/2,3/2,1/2)=0\uparrow\varpi_4=(13/2,5/2,3/2,1/2)
% \uparrow\varpi_1=(13/2,7/2,3/2,1/2)\uparrow\varpi_3=(7,3,2,1)
% \uparrow(15/2,5/2,3/2,1/2)=2\varpi_4\uparrow
% \varpi_1+\varpi_4=(15/2,7/2,3/2,1/2)\uparrow\varpi_2=(15/2,7/2,5/2,1/2)$

Remains the case $p>2$.
To see that $\nabla_G(\lambda)\to\nabla_H(\res^{B}_{B\cap H}\lambda)$
is surjective for $\lambda=\sigma_1$, it suffices to do this for
$\lambda=\varpi_1$ and $\lambda=\varpi_3$.
For $p>3$ one could now use that
$\nabla_H(\res^{B}_{B\cap H}\lambda)$
is irreducible for both $\lambda=\varpi_1$ and $\lambda=\varpi_3$,
because each of the dominant weights of $\nabla_H(\res^{B}_{B\cap H}\lambda)$
is in a different orbit
under the affine Weyl group.
% The dominant weights one has to check here are $0$, $\varpi_4$, $\varpi_1$,
% $\varpi_3$
% in the notations of Bourbaki for {\Ff}.
% They are in the closure of the bottom alcove for $p\geq13$.
% For $p=5$, 7, 11 one just transforms them back to
% the closure of the bottom alcove.
% It took us 33 reflections in total.

But we need an argument that works for $p\geq3$.
Now $\nabla_G(\varpi_1)$ is a miniscule representation of dimension 27,
and $\nabla_H(\varpi_4)=\nabla_H(\res^B_{B\cap H}\varpi_1)$ has dimension 26.
There are 24 short roots and they have multiplicity one in
$\nabla_H(\varpi_4)$.
So the map from $M:=\nabla_G(\varpi_1)$ to $\nabla_H(\varpi_4)$ hits at least
24 dimensions and its kernel consists of $H$-invariants.
Indeed there are three weights of $\nabla_G(\varpi_1)$ that restrict to zero.
In Bourbaki notation
they are $\zeta_1=1/6(\epsilon_8-\epsilon_7-\epsilon_6)+
1/2(-\epsilon_1+\epsilon_2+\epsilon_3-\epsilon_4-\epsilon_5)$,
$\zeta_2=1/6(\epsilon_8-\epsilon_7-\epsilon_6)+
1/2(\epsilon_1-\epsilon_2-\epsilon_3+\epsilon_4-\epsilon_5)$,
$\zeta_3=-1/3(\epsilon_8-\epsilon_7-\epsilon_6)+\epsilon_5$.
Put $\zeta_4=1/6(\epsilon_8-\epsilon_7-\epsilon_6)+
1/2(\epsilon_1-\epsilon_2-\epsilon_3-\epsilon_4+\epsilon_5)$,
$\zeta_5=1/6(\epsilon_8-\epsilon_7-\epsilon_6)+
1/2(-\epsilon_1+\epsilon_2-\epsilon_3+\epsilon_4-\epsilon_5)$.
Then $X_{\alpha_1}$ induces an isomorphism $M_{\zeta_3}\to M_{\zeta_4}$
and it annihilates $M_{\zeta_1}+M_{\zeta_2}$.
Similarly $X_{\alpha_6}$ induces an isomorphism $M_{\zeta_2}\to M_{\zeta_4}$
and annihilates $M_{\zeta_1}+M_{\zeta_3}$.
The same space is annihilated by $X_{\alpha_3}$, which induces an isomorphism
$M_{\zeta_2}\to M_{\zeta_5}$.
Finally
$X_{\alpha_5}$ induces an isomorphism $M_{\zeta_1}\to M_{\zeta_5}$
and annihilates $M_{\zeta_2}+M_{\zeta_3}$.

It follows that in $M_{\zeta_1}+M_{\zeta_2}+M_{\zeta_3}$ there is just
a one dimensional subspace of vectors annihilated by both
$X_{\alpha_1}+X_{\alpha_6}$ and $X_{\alpha_3}+X_{\alpha_5}$.
(These two operators come from the Lie algebra of H.)
We conclude that $\res_H^GM$ has a good filtration and that
$M\to\nabla_H(\varpi_4)$
is surjective. As $p>2$, we then also have that $M\wedge M$ and
$\res_H^G(M\wedge M)$ have a good filtration. It then
follows from the character that
$M\wedge M=\nabla_G(\varpi_3)$.
(We use the program \texttt{LiE} \cite{LiE}.)
% \begin{verbatim}
% > alt_tensor(2,X[1,0,0,0,0,0],E6)
%      1X[0,0,1,0,0,0]
% \end{verbatim}
So $\res_H^G\nabla_G(\varpi_3)$ has a good filtration and therefore maps onto
$\nabla_G(\res_{B\cap H}^B\varpi_3)$.

Summing up, we have shown
\begin{Theorem}
The {\Es}-{\Ff} pair is a Donkin pair.
In fact it satisfies the
pairing criterion.
\end{Theorem}

\begin{Remark}
When Steve Donkin received this proof, he proceeded to show that one
could also prove {\Es}-{\Ff} to be a Donkin pair with the `ancient methods'
of his book \cite{Donkin}. Of course he had to treat more
representations than we do.
We will use his method in the last section
to treat the remaining
cases of Brundan's conjecture in characteristic $p>2$, where we have no
alternative yet.
\end{Remark}

\section{Induction and canonical splitting}
We finish the discussion of canonical splittings
with an analogue of proposition \cite[5.5]{Mathieu T}.
It makes a principle from \cite{Mathieu G} more explicit.
The result was explained to us by O. Mathieu at a reception of the mayor
of Aarhus in August 1998.
It shows once more that canonical splittings combine well with Demazure
desingularisation of Schubert varieties.

\begin{Proposition}\label{inducedcanonical}
Let $X$ be a projective $B$-variety with canonical splitting.
Let $P$ be a minimal parabolic. Then $P\times^BX$ has a canonical splitting.
\end{Proposition}
\begin{Corollary}
The same conclusion holds for any parabolic subgroup.
\end{Corollary}
\paragraph{Proof}If $P$ is not minimal,
take a Demazure resolution
$Z=P_{1} \times^B P_{{2}} \times^B \cdots \times^B P_{r}/B$ of $P/B$
and apply the proposition to
get a canonical splitting on
$P_{1} \times^B P_{{2}} \times^B \cdots \times^B P_{r} \times^BX$.
Then push the splitting forward (\cite[Prop.~4]{Mehta-Ramanathan})
to $P\times^BX$.\qed
\paragraph{Proof of Proposition}
We use notations as in \cite[Ch.~4, A.4]{vdK}.
Let $\zeta$ be the highest weight of $\St$
and $s$ the simple reflection corresponding with $P$.
One checks as in \cite[A.4.6]{vdK} that
$$\sEndF(P\times^BX)=(P\times^B\sEndF(X))\otimes\pi^*\cl(s\zeta-\zeta),$$
where
$\pi:P\times^BX\to P/B$.
We are given a map $\phi:k_\zeta\otimes\St\to\EndF(X)$.
The required map $\psi:k_\zeta\otimes\St\to \EndF(P\times^BX)$
may be constructed by composing maps
\begin{eqnarray*}
k_\zeta\otimes\St&\cong&\\
k_{-s\zeta}\otimes\ind_B^P(k_{\zeta+s\zeta}\otimes\St)&\to&\\
k_{-s\zeta}\otimes\ind_B^P(k_{s\zeta}\otimes\EndF(X))&\cong&\\
k_{-s\zeta}\otimes
H^0(P\times^BX,P\times^B(\sEndF(X)[s\zeta]))&\cong&\\
\EndF(P\times^BX,B\times^BX)&\to&\\
\EndF(P\times^BX)&&
\end{eqnarray*}
Here $k_{-s\zeta}$ is identified with the weight space of weight
$-s\zeta$ of $$H^0(P\times^BX,\pi^*\cl(-\zeta)).$$
An element of that weight space has divisor
$(p-1)B\times^BX=(p-1)X$.

To see that the image of $\psi$ is not in the kernel of
$$\epsilon_{P\times^BX}:\EndF(P\times^BX)\to k,$$
it suffices to show that the diagram
$$ \begin{array}{ccccc}
 k_\zeta\otimes\St & {\to} &
\EndF(P\times^BX,B\times^BX)  & {\to} &
\EndF(P\times^BX)\\[.2em]
 {\|} & & {\downarrow} & & {\downarrow}\\
 k_\zeta\otimes\St  & \stackrel{\phi}{\to} &
\EndF(X) & {\to} & k
\end{array}
$$
commutes.
Now
$$ \begin{array}{ccc}
 k_{-s\zeta}\otimes\ind_B^P(k_{\zeta+s\zeta}\otimes\St) & {\longrightarrow}
&
k_{\zeta}\otimes\St \\[.2em]
 {\downarrow} & & {\downarrow}\\
k_{-s\zeta}\otimes\ind_B^P(k_{s\zeta}\otimes\EndF(X)) & {\longrightarrow} &
\EndF(X)
\end{array}
$$
commutes and by restricting to the trivial fibration
$BsB\times^BX\to BsB/B$ one shows through the following lemma
that the bottom map in this last diagram
agrees with the map that factors through $\EndF(P\times^BX,B\times^BX)$.\qed

\begin{Lemma}
Let $A$ be a commutative $k$-algebra. Then
$$\EndF(A[t])=\EndF(A)\otimes\EndF(k[t])=\EndF(A)\otimes k[t]$$ and
the map $\EndF(A[t],(t))\to\EndF(A)$ is induced by the map
$$t^{p-1}k[t]=t^{p-1}*\EndF(k[t])=\EndF(k[t],(t))\to\EndF(k)=k$$
which sends $t^{p-1}f(t)$ to $f(0)$.
\end{Lemma}

\paragraph{Proof}
Straightforward, provided one keeps in mind how $\EndF(R)$
is an $R$-module
(\cite[4.3.3]{vdK}).
Compare also \cite[A.4.5]{vdK}.\qed

\section{More Donkin pairs}
In this section we do not use the pairing criterion.
Instead we return to the methods of Donkin's book \cite{Donkin}, combined
with computer calculations of characters, of linkage,
and of the Jantzen sum formula.

Let $G$, $H$ be as before, with $G$ simply connected.
In fact $H$ will be the commutator subgroup of
the group of fixed points of an involution of $G$
which leaves invariant the maximal torus $T$ and the Borel subgroup $B$.
We refer to \cite{Springer} for the classification of the possibilities,
assuming $p>2$. (Of course involutions of the simply connected $G$ are
lifted \cite[9.16]{Steinberg}
from the involutions of the corresponding adjoint group, which are
treated in \cite{Springer}.)
% To see what the lift does on the torus, one notes that the action on
% coroots is determined by the action on roots, which may be read off
% in the adjoint group.

\begin{Remark}\label{pass to commutator}
Let $H$ be the fixed point group of an involution that leaves $T$ and $B$
invariant in the simply connected semisimple $G$.
Then $H$ is connected reductive by \cite[8.2]{Steinberg}.
Now an $H$-module has good filtration in the sense of \cite{Donkin}
if and only if its
restriction to the commutator subgroup of $H$ has good filtration.
That is why we look only at semisimple subgroups $H$.
\end{Remark}

Let $\cm$ denote the set of finite dimensional $G$-modules $M$
with good filtration
for which $\res^G_HM$ has good filtration.
Let $\cs$ denote the set of dominant weights $\lambda$ of $G$ so that
$\nabla_G(\lambda)\in \cm$.
As always we try to show that all dominant weights of $G$ are in $\cs$.
For this purpose we recall some useful lemmas.

\begin{Lemma}\label{closure}
\begin{enumerate}
\item If $M_1\oplus M_2\in\cm$, then $M_1\in \cm$.
\item If $0\to M'\to M\to M''\to 0$ is exact, and $M'\in \cm$, then
$M\in \cm$ if and only if $M''\in \cm$.
\item If $M_1$, $M_2\in\cm$, then $M_1\otimes M_2\in\cm$.
\end{enumerate}
\end{Lemma}

If $M$ is a $G$-module with good filtration, write $\suppd(M)$ for the set
of dominant weights $\lambda$ so that $\nabla_G(\lambda)$ occurs as a layer
in a good filtration of $M$.
We order the dominant weights of $G$ by the partial order in which
$\mu\leq\lambda$ if and only if $\lambda-\mu$ is in the closed cone
spanned by the positive roots.
In particular, if $\lambda$ is dominant, then $0\leq\lambda$, and
all dominant weights $\mu$ of $\nabla_G(\lambda)$ satisfy $\mu\leq\lambda$.
We say that a filtration of $M$ is a good filtration \emph{adapted} to the
partial order if there are $\lambda_i$ so that
the $i$-th layer is a direct sum of copies of
$\nabla_G(\lambda_i)$, and $i\leq j$ if $\lambda_i\leq \lambda_j$.
(So we still call it a good filtration, even though $\nabla_G(\lambda_i)$
may have multiplicity in the $i$-th layer.)
If $M$ has a good filtration, then it also has one adapted to the partial
order, by the proof of \cite[II 4.16]{Jantzen}.

\begin{Lemma}\label{auxmodule}
Let $M\in\cm$ and $\lambda\in\suppd(M)$.
Assume for every weight $\mu$ in $ \suppd(M)$, distinct from $\lambda$,
that one of the following holds
\begin{enumerate}
\item $\mu<\lambda$ and $\mu\in\cs$.
\item $\mu$ and $\lambda$ are in different orbits under the affine Weyl group.
\end{enumerate}
Then $\lambda\in\cs$.
\end{Lemma}

\paragraph{Proof}We may replace $M$ by an indecomposable direct summand
$M_1$ with $\lambda\in\suppd(M_1)$. The linkage principle tells
that we thus get rid of the second
possibility in the lemma. Then in a good filtration adapted to the
partial order, the module
$\nabla_G(\lambda)$ occurs only as a summand of the top layer, which is in
$\cm$ by lemma \ref{closure}.\qed

\begin{Lemma}\label{irreducible}
Let $\lambda$ be a dominant weight of $G$.
If $\lambda$ is in the bottom alcove, or if the Jantzen sum formula yields
zero, then $\nabla_G(\lambda)$ is irreducible.
\end{Lemma}

\paragraph{Proof}
See \cite[II Cor.~5.6 and 8.21]{Jantzen}\qed\par

%  In practice we will need this with $\lambda$ a sum of at most two,
%  not necessarily distinct, fundamental weights. Now for the types
%  D8, E7, A5, C4, A7, D6, C3, one has
%  $\langle\lambda_i,\alpha_0^\wedge\rangle\leq4$ and
%  $\langle\rho,\alpha_0^\wedge\rangle\leq17$, so we never need to go
%  beyond 29, which is still covered by the applet.

\subsection{The pairs $E_8$, $D_8$ and $E_8$, $E_7A_1$}

Say $G$ is of type $E8$ in characteristic $p>2$
and $H$ is the fixed point group of an involution.
There are two cases, up to conjugacy. One may have $H$ of type $D_8$ or one
may have $H$ of type $E_7A_1$.

In either case we wish to show that $G$, $H$ is a Donkin pair.
In other words, we want that all dominant weights are in $\cs$.
We will argue by induction along the partial order. Thus when trying to
prove that $\lambda\in \cs$, we shall always assume that $\mu\in \cs$ for
$\mu<\lambda$.
Of course the zero weight is in $\cs$, so say $\lambda$ is nonzero.
If $\lambda$ is not a fundamental weight, write $\lambda=\lambda_1+\lambda_2$
where $\lambda_i$ are nonzero dominant weights. As $\lambda_i<\lambda$, we
may apply lemma \ref{auxmodule} with
$M=\nabla_G(\lambda_1)\otimes\nabla_G(\lambda_2)$ to conclude that
$\lambda\in\cs$.

Remain the fundamental weights. Observe that
$\varpi_8<\varpi_1<\varpi_7<\varpi_2<\varpi_6<\varpi_3
<\varpi_5<\varpi_4$. But we will not discuss
them in this exact order.

To see that $\varpi_8\in\cs$ we compute the character of
$\res^G_H\nabla_G(\varpi_8)$
with the program \texttt{LiE} \cite{LiE},
decompose this character in terms of Weyl characters, and use lemma
\ref{irreducible} to see that $\res^G_H\nabla_G(\varpi_8)$
has a composition series whose factors are
irreducible (dual) Weyl modules. Here we use a Java applet of Lauritzen
for the Jantzen sum formula.

To see that $\varpi_1\in\cs$ we may argue the same way if $H$ is of type $D_8$.
If $H$ is of type $E_7A_1$, let $K$ be the subgroup of type $E_7$ in $H$,
and $F$ the subgroup of type $A_1$.
Then $G$, $K$ is a Donkin pair (Levi subgroup case), so we may consider
a good filtration adapted to the partial order (on weights for $K$)
of $\res_K^G\nabla_G(\varpi_1)$.
It is a filtration by $H$-modules. It suffices to show that its layers
have good filtration as $H$-modules. Let $N$ be such a layer.
Its character is the character of some $\nabla_H(\lambda_1,\lambda_2)$
where $\lambda_1$ is a dominant weight for $K$ and $\lambda_2$ is one
for $F$. (From now on we do not mention the
computer calculations that are needed
to support such statements.)
Moreover, $\nabla_F(\lambda_2)$ is irreducible, so that the
$K\cap B$-socle of $N$ is an irreducible $F$-module. It follows that the
natural map $N\to \nabla_H(\lambda_1,\lambda_2)$ is an isomorphism, and
thus $\varpi_1\in\cs$.

To see that $\varpi_2\in\cs$, we apply lemma \ref{auxmodule} with
$M=\nabla_G(\varpi_1)\otimes\nabla_G(\varpi_8)$. (Note $\varpi_1$,
$\varpi_8<\varpi_2$, so that indeed $M\in\cm$ by the
inductive assumption.)
One may find the necessary statement about linkage in Donkin's book. (This is
no accident, as we follow him in our choices.) We also checked the
non-linkage with
a straightforward Mathematica program.

To see that $\varpi_7\in\cs$, we similarly use $M=\nabla_G(\varpi_8)\wedge
\nabla_G(\varpi_8)$. (Recall $p>2$.)
To get $\varpi_3\in\cs$,
use $M=\nabla_G(\varpi_1)\wedge
\nabla_G(\varpi_1)$.
To get $\varpi_4\in\cs$,
use $M=\nabla_G(\varpi_2)\wedge
\nabla_G(\varpi_2)$.
To get $\varpi_5\in\cs$,
use $M=\nabla_G(\varpi_1)\otimes
\nabla_G(\varpi_2)$ if $p=3$, and $M=\wedge^4\nabla_G(\varpi_8)$ if $p>3$.
To get $\varpi_6\in\cs$,
use $M=\nabla_G(\varpi_1)\otimes
\nabla_G(\varpi_1)$ if $p=3$, and $M=\wedge^3\nabla_G(\varpi_8)$ if $p>3$.

\subsection{The pair  $E_6$,  $A_5A_1$}
Let $G$ be the simply connected group of type $E_6$ in characteristic $p>2$
and let $H$ be the fixed point group of such an inner involution that $H$
is of type $A_5A_1$ and the involution commutes with the graph automorphism.
We wish to show again this is a Donkin pair.
We argue as in the $E_8$, $E_7A_1$ case.

If $\lambda=\varpi_1$ or $\varpi_2$, then we argue with socles as we did
to show $\varpi_1\in\cs$ for the $E_8$, $E_7A_1$ pair.

To treat $\varpi_3$ we use $M=\nabla_G(\varpi_1)\wedge
\nabla_G(\varpi_1)$.
To get $\varpi_4\in\cs$,
use $M=\nabla_G(\varpi_2)\wedge
\nabla_G(\varpi_2)$.
The remaining two fundamental weights are in $\cs$ by symmetry.

\subsection{The pair  $E_6$,  $C_4$}
Let $G$ be the simply connected group of type $E_6$ in characteristic $p>2$
and let $H$ be the fixed point group of such an outer involution that $H$
is of type $C_4$ and the involution commutes with the graph automorphism.

The module $\nabla_G(\varpi_2)$ is the Lie algebra of the adjoint form of $G$.
Its restriction $\res^G_H\nabla_G(\varpi_2)$
has a six dimensional weight space for the weight
zero, just like $\nabla_G(\varpi_2)$ itself.
We claim that it contains no nonzero invariant. Indeed we may choose the
involution so that
$X_{\alpha_3}+X_{\alpha_5}$,
$X_{\alpha_1}+X_{\alpha_6}$,
$X_{\alpha_3+\alpha_4}+X_{\alpha_5+\alpha_4}$, $X_{\alpha_2}$
are in the Lie algebra of $H$,
where we have put $X_{\alpha_3+\alpha_4}=[X_{\alpha_3},X_{\alpha_4}]$
and $X_{\alpha_5+\alpha_4}=[X_{\alpha_4},X_{\alpha_5}]$.
The only element in the weight zero weight space of $\nabla_G(\varpi_2)$
that is annihilated by all these elements is the zero vector.
Now $\res^G_H\nabla_G(\varpi_2)$ contains an irreducible $\nabla_H(2\varpi_1)$
and the quotient by that submodule is either irreducible, or $p=3$ and
there are two composition factors, one of which is one dimensional.
(This also uses the Jantzen sum formula.) As there is no invariant
in $\res^G_H\nabla_G(\varpi_2)$ and there is no extension between
$\nabla_H(2\varpi_1)$ and the other composition factors
(\cite[II 4.13, 4.14]{Jantzen}), we
get $\varpi_2\in\cs$.

As $\res^G_H\nabla_G(\varpi_1)$ is irreducible, we also have $\varpi_1\in\cs$.
The rest goes as for the pair  $E_6$,  $A_5A_1$.

\subsection{The pairs $E_7$, $A_7$ and $E_7$, $D_6A_1$}

Say $G$ is simply connected of type $E7$ in characteristic $p>2$
and $H$ is the fixed point group of an involution.
There are two cases, up to conjugacy. One may have $H$ of type $A_7$ or one
may have $H$ of type $D_6A_1$.

We argue as before. If $H$ is of type  $D_6A_1$ we show that
$\varpi_1,\varpi_2,\varpi_7\in\cs$ by the argument with socles used to show
$\varpi_1\in\cs$ for the $E_8$, $E_7A_1$ pair. If $H$ is of type $A_7$
we have $\varpi_1,\varpi_7\in\cs$ for the same reason,
involving the sum formula, as why
$\varpi_8\in\cs$ for the $E_8$, $D_8$ pair.
If $p=7$ we see in the same manner that $\varpi_2\in\cs$.
If $p\neq7$ use $M=\nabla_G(\varpi_1)\otimes\nabla_G(\varpi_7)$
to get $\varpi_2\in\cs$.
To get $\varpi_3\in\cs$ use $M=\nabla_G(\varpi_1)\wedge\nabla_G(\varpi_1)$.
To get $\varpi_4\in\cs$ use $M=\nabla_G(\varpi_2)\wedge\nabla_G(\varpi_2)$.
To get $\varpi_5\in\cs$ use  $M=\wedge^3\nabla_G(\varpi_7)$ if $p\neq3$
and $M=\nabla_G(\varpi_1)\wedge\nabla_G(\varpi_2)$ otherwise.
To get $\varpi_6\in\cs$ use $M=\nabla_G(\varpi_7)\wedge\nabla_G(\varpi_7)$.

\subsection{The pairs $F_4$, $B_4$ and $F_4$, $C_3A_1$}
Say $G$ is of type $F_4$ in characteristic $p>2$
and $H$ is the fixed point group of an involution.
There are two cases, up to conjugacy. One may have $H$ of type $B_4$ or one
may have $H$ of type $C_3A_1$.

We argue as before. If $H$ is of type  $C_3A_1$ we show that
$\varpi_1,\varpi_4\in\cs$ by the argument with socles used to show
$\varpi_1\in\cs$ for the $E_8$, $E_7A_1$ pair.
If $H$ is of type $B_4$ we have $\varpi_1,\varpi_4\in\cs$ for
the same reason, involving the sum formula, as why
$\varpi_8\in\cs$ for the $E_8$, $D_8$ pair.
To get $\varpi_3\in\cs$ use $M=\nabla_G(\varpi_4)\wedge\nabla_G(\varpi_4)$.
To get $\varpi_2\in\cs$ use $M=\nabla_G(\varpi_1)\wedge\nabla_G(\varpi_1)$.

\begin{Theorem}[Brundan's Conjecture \cite{Brundan}]
Let $G$ be semisimple simply connected.
If either
\begin{enumerate}
\item[(i)] $H$ is the centralizer of a graph automorphism of $G$; or
\item[(ii)] $H$ is the centralizer of an involution of $G$ and the
characteristic is at least three,
\end{enumerate}
then $G$, $H$ is a Donkin pair.
\end{Theorem}

\subsubsection*{Proof}
We have either a Levi subgroup case, first settled
in \cite{Donkin} (see also remarks \ref{levi},
\ref{pass to commutator}), or a case treated in \cite{Brundan}, or a case
treated
above, up to conjugacy.\qed

\begin{Remark}
Of course
we would much prefer a case-free proof, based on the pairing criterion say.
\end{Remark}

\end{document}